\documentclass[journal]{IEEEtran}

\usepackage{bm}
\usepackage{wrapfig}
\usepackage{float}
\usepackage{amsmath}
\usepackage{amssymb}
\usepackage{multirow}
\usepackage{amsfonts}
\usepackage{mathrsfs}
\usepackage{graphicx}
\usepackage{epstopdf}
\usepackage{subcaption}
\usepackage{adjustbox}
\usepackage{textcomp}
\usepackage{autobreak}
\usepackage{tabularx}
\usepackage{comment}
\usepackage{caption}
\usepackage{lettrine}
\usepackage{cite}
\usepackage{setspace}
\usepackage{amsmath}
\usepackage{titlesec}
\usepackage{lipsum}
\usepackage{glossaries}
\usepackage[utf8]{inputenc}
\usepackage{tikz}
\usepackage{import}
\usepackage{dirtytalk}
\usepackage{algorithm}
\usepackage{algpseudocode}

\titlespacing*{\section}
{0pt}{1.5ex plus 1ex minus .1ex}{2.0ex plus .1ex}

\begin{document}

\title{Cost-Aware Bound Tightening for Constraint Screening in AC OPF}

\author{Mohamed~Awadalla and Fran\c{c}ois~Bouffard,~\IEEEmembership{Senior Member,~IEEE}%
\thanks{This work was supported in part by IVADO, Montreal, QC, the Trottier Institute for Sustainability in Engineering and Design, Montreal, QC, and the Natural Sciences and Engineering Research Council of Canada, Ottawa, ON.}
\thanks{M. Awadalla and F. Bouffard are with the Department of Electrical and Computer Engineering, McGill University, Montreal, QC H3A~0E9, Canada and with the Groupe d'\'{e}tudes et de recherche en analyse des d\'{e}cisions (GERAD), Montreal, QC  H3T~1J4, Canada (emails: mohamed.awadalla@mail.mcgill.ca; francois.bouffard@mcgill.ca).}
}

\maketitle

\begin{abstract}
The objective of electric power system operators is to determine cost-effective operating points by resolving optimization problems that include physical and engineering constraints. As empirical evidence and operator experience indicate, only a small portion of these constraints are found to be binding during operations. Several optimization-based methods have been developed to screen out redundant constraints in operational planning problems like the optimal power flow (OPF) problem. These elimination procedures primarily focus on the feasible region and ignore the role played by the problem's objective function. This letter addresses the constraint screening problem using the bound tightening technique in the context of the OPF problem formulated with a full ac power flow characterization. Due to the non-convexity of the ac OPF, we investigate line constraint screening under different convex relaxations of the problem, and we evaluate how the economics of the objective function impacts screening outcomes. 
\end{abstract}


\section{Introduction}
\IEEEPARstart{T}{he} optimal power flow (OPF) problem seeks to find optimal operating conditions subject to networks' physical and engineering constraints \cite{Lehmann, Ardakani2}. The OPF incorporates physical constraints to represent power flow phenomena and  constraints such as voltage, generator, and line flow limits. The ac optimal power flow (AC-OPF) problem is both nonlinear, nonconvex and therefore NP-hard \cite{Lehmann}. Furthermore, the number of variables and constraints becomes large as network size increases, posing more computational challenges for large realistic systems.\par

Despite the fact that all constraints in the OPF problem must be satisfied for a solution to be feasible, power system operators' experience and past research have shown only a small proportion of the problem's inequality constraints (especially line flow limits) can be potentially binding (i.e., satisfied with equality) at the optimal solution \cite{Ardakani2}. As a result, it is widely common to screen out redundant constraints in order to reduce the problem size and speed up computations. \par

Umbrella constraint discovery is a pioneer constraint screening algorithm which has been proposed for DC-OPF problems to filter out redundant constraints \cite{Ardakani2}. Another dominant technique is \emph{optimization-based bound tightening} (OBBT) which was proposed primarily to enhance the quality of convex relaxations \cite{Coffrin1,Chen}. The general procedure of OBBT is performed by minimizing or maximizing the desired variable while considering the constraints of the relaxation. This approach has been extended to screen out line flow constraints in DC-OPF, unit commitment \cite{Roald} and AC-OPF \cite{Aquino} problems. Using OBBT, the authors of \cite{Roald, Aquino} identify line flow limits that will never become active by solving one minimization and one maximization optimization problem associated with each line flow limit. In a broader sense, several constraints are satisfied indirectly through other constraints in the problem, allowing them to be eliminated confidently before requesting the solver's assistance.

On the other hand, several convex relaxations of the AC-OPF problem have garnered substantial attention for various reasons. These relaxations include the quadratic convex relaxation (QCR) \cite{Coffrin1}, the semidefinite relaxation (SDR) \cite{Lavaei}, and the second-order cone relaxation (SOCR) \cite{Kocuk}.
Generally, convex relaxations are only approximations; they offer a limit on the best possible global optimal value of the AC-OPF problem. The advantage of convex relaxations with respect to constraint screening problems is that relaxations provide a lower (optimistic) objective bound (assuming the problem is a minimization; the converse applies if one is solving a maximization). Therefore, if the limit provided by a relaxation falls within the set boundaries for the restricted variable, we can be confident that the constraint is unnecessary or redundant in the original problem \cite{Aquino}.

In this letter, we propose a refinement to OBBT for AC-OPF constraint screening. First, instead of attempting to solve OBBT with the nonconvex AC-OPF problem directly, we utilize a convex relaxation to establish an upper limit on the global solution. Second, rather than identifying only redundant line constraints as suggested in \cite{Aquino}, we propose a valid upper bound inequality constraint that embodies prior economical information to filter out more redundant constraints. 
This refinement distinguishes between redundant, active, and inactive constraints in the AC-OPF problem. Finally, we apply the proposed method using state-of-the-art convex relaxations \cite{Bingane1} to compare the relative effectiveness of the our proposal. The main contributions of this letter are twofold:
\begin{enumerate}
\item To investigate constraints screening for the AC-OPF problem and show how different convex relaxations yield varying outcomes. We provide a comparative analysis to determine the conservativeness of various convex relaxations in terms of constraint screening. 
\item To improve the OBBT method to identify not only redundant constraints but also inactive ones, thus leading to lighter AC-OPF formulations.
\end{enumerate}


\section{AC Optimal Power Flow Problem Formulation}

Consider a typical power system which consists of a set of $\mathcal{N}=$ $\{1, \ldots, n\}$ buses (a subset of which have generators denoted by the set $\mathcal{G}$) and $\mathcal{L}=$ $\{1, \ldots, l\}$ branches. Every node $k \in \mathcal{N}$ in the network has three properties, voltage $v_k = V_{k}\angle \delta_k$, power generation $s_{G_k}= p_{G_k}+jq_{G_k}$, and power consumption $s_{D_k}= p_{D_k}+jq_{D_k}$, all of which are complex numbers because of the AC power's oscillatory characteristics. Each branch $\ell \in \mathcal{L}$ has a series admittance $y_{\ell}=g_{\ell}+j b_{\ell}$ and a total shunt susceptance $b_{\ell}^{\prime}$. We denote each bus's shunt conductance and susceptance as $g_k^{\prime}$ and $b_k^{\prime}$, respectively. Each branch $\ell \in \mathcal{L}$ has a \emph{from} end $k$ and a \emph{to} end $m$ such that we denote $\ell=(k, m)$. From these we have active and reactive power flows for each branch $\ell$ leaving its \emph{from} end---$p_{f\ell}, q_{f \ell}$---and corresponding flows at its \emph{to} end---$p_{t \ell}, q_{t \ell}$. Lastly, in the case where a branch $\ell$ is a tap-changing transformer, we denote its off-nominal setting with the symbol $t_\ell$. 
The AC-OPF problem is formulated as specified in \cite{Bingane1}:
\begin{equation}
\label{eq1}
\min \sum_{g \in \mathcal{G}} c_{g 2} p_{G g}^2+c_{g 1} p_{G g}+c_{g 0}
\end{equation}
where variables $\boldsymbol{p}_G, \boldsymbol{q}_G \in \mathbb{R}^{|\mathcal{G}|}, \boldsymbol{p}_f, \boldsymbol{q}_f, \boldsymbol{p}_t, \boldsymbol{q}_t, \boldsymbol{t} \in \mathbb{R}^{|\mathcal{L}|}$, and $\boldsymbol{v} \in \mathbb{C}^{|\mathcal{N}|}$, is subject to:
\begin{align}
\sum_{g \in \mathcal{G}_k} & p_{G g}-p_{D k}  -g_k^{\prime}\left|v_k\right|^2 
\nonumber \\
& =\sum_{\ell=(k,m) \in \mathcal{L}} p_{f \ell}+\sum_{\ell=(m, k) \in \mathcal{L}} p_{t \ell}, \quad \forall k \in \mathcal{N} \label{eq2}
\end{align}
\begin{align}
\sum_{g \in \mathcal{G}_k} & q_{G_g} -q_{D k}+b_k^{\prime}\left|v_k\right|^2 \nonumber \\
& = \sum_{\ell=(k,m) \in \mathcal{L}} q_{f \ell}+\sum_{\ell=(m,k) \in \mathcal{L}} q_{t \ell}, \quad \forall k \in \mathcal{N} \label{eq3}
\end{align}
\begin{align}
\frac{v_k}{t_{\ell}} \left [ \left ( j \frac{b_{\ell}^{\prime}}{2} + y_{\ell} \right ) \right . & \left . \frac{v_{k}}{t_{\ell}} - y_{\ell} v_m \right ]^*  \nonumber \\
  & =  p_{f\ell}+j q_{f\ell}, \quad \forall \ell=(k,m) \in \mathcal{L} \label{eq4} \\
v_m  \Bigg [ -y_{\ell} \frac{v_{\ell}}{t_{\ell}} + \Big ( & j \frac{b_{\ell}}{2}  + y_{\ell} \Big ) v_m \Bigg ]^*  \nonumber \\
 & = p_{t\ell}+ j q_{t\ell}, \quad \forall \ell=(k,m) \in \mathcal{L} \label{eq5}
 \end{align}
 \begin{equation}
 \label{eq6}
\underline{p}_{G_g} \leq p_{G_g} \leq \overline{p}_{G_g}, \; \underline{q}_{G g} \leq q_{G_g} \leq \overline{q}_{G g}, \quad \forall g \in \mathcal{G}
 \end{equation}
 \begin{equation}
 \label{eq7}
 \left|p_{f \ell}+ j q_{f \ell}\right| \leq \bar{s}_{\ell}, \; \left|p_{t \ell}+ j q_{t \ell}\right| \leq \bar{s}_{\ell},  \quad \forall \ell \in \mathcal{L} 
 \end{equation}
 \begin{equation}
 \label{eq9}
\underline{v}_k \leq\left|v_k\right| \leq \overline{v}_k, \quad \forall k \in \mathcal{N}
 \end{equation}
 \begin{equation}
 \label{eq10}
 \angle v_1=0 
 \end{equation}
The objective function \eqref{eq1} minimizes the total production cost, with $c_{g2}$, $c_{g1}$ and $c_{g0}$ denoting the coefficients for a quadratic cost function for all $g \in \mathcal{G}$. Constraints \eqref{eq2}--\eqref{eq3} represent the nodal power balance equations for active and reactive powers, respectively. Constraints \eqref{eq4}--\eqref{eq5} represent the active and reactive power flow equations in each branch. In addition to these physical laws, operational constraints \eqref{eq6}--\eqref{eq9} are required in AC power flows. Constraints \eqref{eq6}--\eqref{eq7} impose limits on generator active and reactive power outputs and line thermal limits, respectively. We assume that $\underline{v}_k>0$ for all $k \in \mathcal{N}$ in \eqref{eq9}. Constraint \eqref{eq10} specifies, without loss of generality, node $k=1$ as the reference.

Constraints \eqref{eq4} and \eqref{eq5} are nonlinear and non-convex; this makes problem \eqref{eq1}--\eqref{eq10} difficult to solve and, in fact, NP-hard \cite{Lehmann}. Applying local methods to this problem provides no guarantees on the optimality of any solution found. Moreover, it is intractable to solve to global optimality for large-scale instances. Hence, techniques aiming at convexifying and reducing the dimensions of this problem are part of one's arsenal in the hope of efficiently finding a good local optimal solution to this problem.


\section{Line Flow Constraint Screening}
This letter proposes a screening method for line thermal limits adapted to the AC-OPF problem. The method uses optimizations to determine each line's minimum and maximum power flow values subject to all other constraints. If those flows are within the specified line flow limits, the limits are considered to be redundant and can be eliminated. However, if the flows reach established limits, the constraint is non-redundant. The method first assesses the binding potential of constraints based on the other problem constraints without considering generation cost functions.

Here we formulate the OBBT problem. It involves the solution of two maximization and two minimization problems for each line flow constraint across a range of load fluctuations. The constraint screening problem objective \eqref{eq11} aims to minimize and to maximize active and reactive power flows while considering all the remaining constraints of the original OPF \eqref{eq2}--\eqref{eq10}, including added load variability constraints \eqref{eq12}--\eqref{eq13}. The optimization problem includes additional decision variables for every active load $p_{D_k}$ and reactive load $q_{D_k}$ for each $k \in \mathcal{N}$. We assume that the power losses are small in typical meshed transmission networks, which means $p_{f \ell} \approx-p_{t \ell}$, while a similar assumption applies to reactive power flows. Hence, we consider only one line end (\emph{to}) for the optimization problems in \eqref{eq11}. The parameter $\Delta$ specifies load uncertainty ranges, and $p_{D_k}^\circ$ and $q_{D_k}^\circ$ refer respectively to nominal nodal active and reactive power loadings. 
\begin{equation}
\min_{\boldsymbol{v},\boldsymbol{t}, \boldsymbol{p}_{G},\boldsymbol{q}_{G},\boldsymbol{p}_{D},\boldsymbol{q}_{D}} / \max_{\boldsymbol{v},\boldsymbol{t}, \boldsymbol{p}_{G},\boldsymbol{q}_{G},\boldsymbol{p}_{D},\boldsymbol{q}_{D}} p_{t\ell} / q_{t\ell}   \\ \label{eq11}
\end{equation}
\\ \text {subject to:}
\begin{equation}
\text {Constraints} \quad \eqref{eq2}-\eqref{eq10} \\ \label{eq14}
\end{equation}
\begin{equation}
\quad(1-\Delta) p_{D_k}^o \leq p_{D_k} \leq(1+\Delta) p_{D_k}^o, \quad \forall k \in \mathcal{N} \\ \label{eq12}
\end{equation}
\begin{equation}
(1-\Delta) q_{D_k}^o \leq q_{D_k} \leq(1+\Delta) q_{D_k}^o, \quad \forall k \in \mathcal{N} \\ \label{eq13}
\end{equation}

When the optimal solution of \eqref{eq10}--\eqref{eq13} (whether one is maximizing or minimizing and is optimizing either one of $p_{t\ell}$ or $q_{t\ell}$) is such that $|p_{t\ell}^\star + jq_{t\ell}^\star| = \bar{s}_\ell$, it indicates that the flow limit for line $\ell$ is non-redundant. On the other hand, if flows are strictly within their allowable range, the corresponding bounds are redundant and could be ignored when solving \eqref{eq1}--\eqref{eq10}. 


\begin{table*}[t!]
\captionsetup{font=scriptsize}
\caption{\sc{Proportion (\%) of Redundant Line Flow Constraints Identified for Different Convex Relaxations}} 
\centering 
\scalebox{0.75}{
\begin{tabular}{c|rrrr|rrrr|rrrr} 
\hline
Method & SDR & TCR (\%) & QCR (\%) & SOCR (\%) & SDR & TCR (\%) & QCR (\%) & SOCR (\%) & SDR & TCR & QCR & SOCR  \\ [1.0ex] 
\hline 
Test case           & \multicolumn{4}{c}{Without Bound Tightening (WTB)}  & \multicolumn{4}{c}{With Bound Tightening (WB)} & \multicolumn{4}{c}{ (WB-WTB)/WTB (\%)}\\ [1.0ex] 
\hline  
case5\_pjm          & 50     & 16 ($-$67\%)     & 0 ($-$100\%)        & 0 ($-$100\%)              & 50     & 33 ($-$33\%)  & 0 ($-$100\%)     & 0 ($-$100\%)           & 0      & 100  & 0   & 0           \\[1ex]
case14\_ieee        & 100    & 100 ($-$0\%)        & 95 ($-$5\%)         & 75 ($-$25\%)              & 100    & 100 ($-$0\%)        & 95 ($-$5\%)      & 90 ($-$10\%)           & 0      & 0    & 0   & 20          \\[1ex]
case24\_ rts        & 90  & 66 ($-$27\%)  & 68 ($-$24\%)  & 21 ($-$77\%)        & 95  & 95 ($-$0\%)      & 87 ($-$8\%)   & 45 ($-$53\%)  & 6   & 44    & 27 & 113   \\[1ex]
case39\_epri        & 74  & 48 ($-$35\%) & 41 ($-$44\%)  & 30 ($-$59)       & 91   & 61 (33\%)   & 41 (55\%)     & 35 ($-$62)    & 24  & 27  & 0       & 14   \\[1ex]
case57\_ieee        & 99   & 95 ($-$4\%)  & 86 ($-$13\%)  & 69 ($-$35\%)       & 100    & 98 ($-$3\%)     & 91 ($-$9\%)    & 83 ($-$18\%)   & 1   & 3  & 6     & 20      \\[1ex]
case73\_ieee\_rts   & 75      & 33 ($-$56\%)  & 56 ($-$26\%)  & 18 ($-$77\%)   & 92  & 67 ($-$27\%)  & 69 ($-$25\%)   & 28 ($-$70\%)     & 22  & 100  & 24      & 57   \\[1ex]
case118\_ieee       & 52   & 19 ($-$64\%)  & 23 ($-$55\%)  & 13 ($-$57\%)      & 56  & 20 ($-$63\%)   & 26 ($-$54\%)    & 13 ($-$77\%)  & 8   & 9  & 12   & 0       \\[1ex]
\hline
\textbf{Average}    & $\mathbf{77}$   & $\mathbf{54 \;(-36\%)}$   & $\mathbf{53 \; (-38\%)}$   & $\mathbf{32 \; (-63\%)}$  & $\mathbf{83}$  & $\mathbf{68 \; (-23\%)}$  & $\mathbf{59 \; (-36\%)}$   & $\mathbf{42 (-56\%)}$  & $\mathbf{-9}$   & $\mathbf{-40}$   & $\mathbf{-10}$     & $\mathbf{-32}$   \\[1ex]
\hline
\end{tabular}}
\label{table1_letter} 
\end{table*}

\section{Cost-Driven Constraint Screening}

To further the power of OBBT in screening line flow constraints, we propose to add an extra valid inequality whose role is to capture the effect of the original problem's objective function \eqref{eq1} as part of constraint screening. For this, we add the constraint
\begin{equation}
\label{eq16}
\sum_{g \in \mathcal{G}} c_{g 2} p_{G g}^2+c_{g 1} p_{G g}+c_{g 0} \leq \bar{C}
\end{equation}
to \eqref{eq11}--\eqref{eq13}. This constraint, by putting an upper bound on the operational cost, ends up identifying line flow constraints which are not only non-redundant but also potentially binding in the original problem. This is the case because \eqref{eq16} limits the allowable power generation in a way similar the objective function \eqref{eq1} is attempting to minimize cost. We note here that the right-hand side of \eqref{eq16} $\bar{C}$ is a function of the forecasted net load in the power system, and it reflects expected system-level operating costs for the given net load forecast. Also, $\bar{C}$ would be set slightly higher than historically-observed system costs for similar loading levels to lower the risk of infeasibility of OBBT.


\section{Results}
This section demonstrates the methods described in Sections III and IV for various PGLib-OPF test cases \cite{2019PGLib}. The implementations of each convex relaxation are based on the CONICOPF package \cite{Bingane1}. We solved all the relaxations in MATLAB using CVX 2.2 with the solver MOSEK 9.1.9 and default precision. All computations were carried out on an Intel Core i7-11700K CPU @ 3.60 GHz and 64~GB of RAM. 

Also, the original AC-OPF problem was solved using MATPOWER-solver MIPS 7.0 to calculate $\bar{C}$. Moreover, to avoid an over-tight cost upper bound, which could lead \eqref{eq11}--\eqref{eq16} to infeasibility, we multiply historical cost values by a conservativeness factor of $102\%$. In all cases here, we consider load variability of $\Delta = 0$ around the nominal load values. All tap changers $\boldsymbol{t}$ are set to their nominal levels and are not optimized.


Finally, the results in Table~\ref{table1_letter} summarize how the proposed constraint screening approach works with (WB) and without (WTB) the addition of the cost bound \eqref{eq16} for the four convex relaxations. Firstly, all convex relaxations show a considerable enhanced constraint screening performance with cost-driven tightening in comparison to the approach without cost-driven tightening. Secondly, the SDR reveals the highest constraint screening ability with an average of 83\% and 77\% superfluous constraint elimination with and without bound tightening, respectively. This outcome complies with the fact that the SDR is the tightest relaxation among those tested. Screening results for all other relaxations are compared to those obtained with SDR. SOCR is the weakest relaxation among the others considering the screening approaches. TCR dominates QCR in five out of seven instances for the WB approach, while for the WTB approach, TCR outperforms QCR in four out of seven instances. In fact, TCR has achieved an optimum enhancement in terms of redundancy removal considering the WB approach in comparison to the WTB approach by an average of 40\%.

\section{Conclusion}
Using the PGLib-OPF test cases with up to 118 buses, we
show that on average, the SDR relaxation yields the best redundant line flow constraint removal performance. Furthermore, for
the test cases with and without bound tightening, TCR is stronger
than QCR. Overall, the proposed approach is set to enhance the solution process of the AC-OPF problem by keeping only its potential constraints which are both non-redundant and active. These improvements in constraint screening performance are achievable thanks to the addition of a cost-driven valid inquality to OBBT.

\bibliography{Reference} 
\bibliographystyle{ieeetr}

\end{document}